\magnification=\magstep1
\hsize=16.5 true cm 
\vsize=23.6 true cm
\font\bff=cmbx10 scaled \magstep1

\font\bffg=cmbx10 scaled \magstep3
\font\bffgg=cmbx10 scaled \magstep4

\parindent0cm
\font\boldmas=msbm10           
\def\Bbb#1{\hbox{\boldmas #1}} 
%%%%%%%%%%%%%%%%%%%%%%%%%%%%%%%%%%%%%%
%%%%%%%%%%%%%%%%%%%%%%%%%%%%%%%%
%%%%ABKUERZUNGEN im Textsatz %%%%
           %
            %
              %
\def\mp{\medskip}              %
            %
               %
%%%%%%%%%%%%%%%%%%%%%%%%%%%%%%%%
%%%ZAHLENMENGEN%%%%%%%%%%%%%%%%%
\font\boldmas=msbm10           %
\def\Bbb#1{\hbox{\boldmas #1}} %
\def\R{\Bbb R}                 %
\def\K{\Bbb K}                 %
%%%%%%%%%%%%%%%%%%%%%%%%%%%%%%%%
 
%%%%%%%%%%%%%%%%%%%%%%%%%%%%%%%%%%%%%%%
\centerline{\bffgg Cantor sets and fields of reals}
\bigskip\medskip
\centerline{\bffg Gerald Kuba}
\bigskip\bigskip\medskip
{\bff 1. Introduction} 
\medskip
The cardinal number (the {\it size}) of a set $\,S\,$
is denoted by $\,|S|\,$. So 
$\,|{\Bbb R}|=|{\Bbb D}|={\bf c}=2^{\aleph_0}\,$ where 
$\;{\Bbb D}\,=\,
\{\,\sum_{n=1}^\infty a_n\,3^{-n}\;\,|\,\;a_n\in\{0,2\}\,\}\;$
is the {\it Cantor ternary set}.
\medskip
As usual, a nonempty set $\,C\subset{\Bbb R}\,$ is a 
{\it Cantor set} if and only if $\,C\,$ is compact and does not 
contain nondegenerate intervals or isolated points.
Equivalently, $\,C\,$ is a 
{\it Cantor set} if and only if there is a continuous 
bijection from $\,{\Bbb D}\,$ onto $\,C\,$.
\medskip
A fundamental question in descriptive analysis of the reals is whether 
a set $\,X\subset{\Bbb R}\,$ contains a Cantor set.
An answer to this question can be important for problems 
concerning size and measure. For example, 
if $\,X\subset{\Bbb R}\,$ contains a Cantor set then (due to $\,|{\Bbb D}|={\bf c}\,$)
the sets $\,X\,$ and $\,{\Bbb R}\,$ are equipollent.
This leads to the important observation that, while  
$\;\aleph_0<|X|<{\bf c}\;$ cannot be ruled out for arbitrary 
sets $\,X\,$, we can be sure that $\;\aleph_0<|A|<{\bf c}\;$
is impossible for {\it closed} subsets $\,A\,$ of the real line.
(Because it is well-known that every uncountable closed subset of $\,{\Bbb R}\,$
must contain a Cantor set.)
Another important example 
is the following. As usual, $\,B\subset{\Bbb R}\,$ is a {\it Bernstein set}
when neither $\,B\,$ nor $\,{\Bbb R}\setminus B\,$ contains a Cantor set.
(Equivalently,                                                           
neither $\,B\,$ nor $\,{\Bbb R}\setminus B\,$ contains an uncountable closed
set.) It is well-known that a Bernstein set is never Lebesgue measurable,
see~[1] Theorem 6.3.8.
\medskip
We are interested in the existence of 
Cantor subsets from a specific algebraic point of view.
Let us call a proper subfield of $\,{\Bbb R}\,$ a {\it Cantor field}
if and only if it contains some Cantor set.
(Notice that a proper subfield of $\,{\Bbb R}\,$
cannot contain the Cantor ternary set $\,{\Bbb D}\,$
because it is well-known that 
$\;\{\,x+y\;|\;x,y\in{\Bbb D}\,\}\,=\,[0,2]\,$.) 
Furthermore, a {\it Bernstein field} is a subfield of $\,{\Bbb R}\,$ 
which is a Bernstein set. Trivially, each Bernstein field is a proper
subfield of $\,{\Bbb R}\,$ and it cannot be a Cantor field.
A trivial consequence of $\,|{\Bbb D}|={\bf c}\,$ is that
$\,|K|={\bf c}\,$ for every Cantor field $\,K\,$. 
We also have $\,|K|={\bf c}\,$ for every Bernstein field $\,K\,$.
Actually, the following is true.
\medskip
(1.1)\qquad{\it $\,|B|={\bf c}\,$ for every Bernstein set $\,B\subset{\Bbb R}\,$.} 
\medskip
There are several ways to verify (1.1). For example, (1.1) is an 
immediate consequence of (2.1) below.
\medskip
For abbreviation let us call fields (taken from some collection) 
{\it incomparable} if no field is isomorphic to 
a subfield of another field. In particular,  
incomparable fields are mutually non-isomorphic.
Our first main result is the following theorem. 
Notice that the field 
$\,{\Bbb R}\,$ is an {\it algebraic extension}
of a subfield $\,K\,$ if $\,\overline{K}={\Bbb C}\,$.
(As usual, $\,\overline{K}\,$ denotes the 
{\it algebraic closure} of the field $\,K\,$.)
\medskip\smallskip
{\bf Theorem 1.} {\it There exist
four families $\,{\cal C}_1,{\cal C}_2,{\cal B}_1,{\cal B}_2\,$
of subfields $\,K\,$ of $\,{\Bbb R}\,$ with $\,\overline{K}={\Bbb C}\,$
such that $\,|{\cal C}_1|=|{\cal C}_2|=|{\cal B}_1|=|{\cal B}_2|=
2^{\bf c}\,$ and all fields in $\,{\cal C}_1\cup{\cal C}_2\,$ 
are {\it Cantor fields} and all fields in $\,{\cal B}_1\cup{\cal B}_2\,$ 
are {\it Bernstein fields} and 
all fields in $\,{\cal C}_1\cup{\cal B}_1\,$ are isomorphic, 
whereas the fields in $\,{\cal C}_2\cup{\cal B}_2\,$ 
are  incomparable.}
\vfill\eject
\bigskip
{\bff 2. Preparation of the proof}
\medskip
For $\,a\in{\Bbb R}\,$ and $\,S\subset{\Bbb R}\,$ 
we write $\;a+S\,:=\,\{\,a+x\;|\;x\in S\,\}\,$.
The following observation implies (1.1) and is also important 
for the proof of Theorem 1.
\medskip
(2.1)\quad {\it Every Cantor set can be partitioned into $\,{\bf c}\,$
Cantor sets.}
\medskip
To verify (2.1) let $\,C\,$ be a Cantor set and $\,f\,$ be 
homeomorphism from $\,{\Bbb D}\,$ onto $\,C\,$.
Now consider the subset
$\;D\,:=\,
\{\,\sum_{n=1}^\infty a_n\,3^{-2n}\;\,|\,\;a_n\in\{0,2\}\,\}\;$
of $\,{\Bbb D}\,$ which is clearly a Cantor set.
Consequently, $\;{\cal P}\,=\,\{\,x+3\cdot D\;|\;x\in D\,\}\;$
is a partition of $\,{\Bbb D}\,$ consisting of 
Cantor sets and $\,|{\cal P}|=|D|={\bf c}\,$.
Thus $\;\{\,f(P)\;|\;P\in{\cal P}\,\}\;$ is a partition of $\,C\,$
as desired, {\it q.e.d.}
\medskip\smallskip
As a consequence of the following lemma, a proper subfield $\,K\,$
of $\,{\Bbb R}\,$ is a Bernstein field if 
$\;K\cap C\not=\emptyset\;$ for every Cantor set $\,C\,$ or, equivalently,
if $\;{\Bbb R}\setminus K\;$ does not contain a Cantor set.
(It is not necessary to check that $\,K\,$ does not contain a Cantor set.)
\medskip
{\bf Lemma 1.} {\it If $\,K\,$ is a proper subfield of $\,{\Bbb R}\,$
which contains a Cantor set then $\,{\Bbb R}\setminus K\,$ contains
a Cantor set as well.}                                      
\smallskip 
{\it Proof.} Let $\,K\,$ be a subfield of $\,{\Bbb R}\,$
and $\,\xi\in{\Bbb R}\setminus K\,$. Then $\;\xi+k_1=k_2\;$ is 
impossible for all $\;k_1,k_2\in K\;$ and hence the translate
$\,\xi+K\,$ is disjoint from $\,K\,$.
Consequently, if $\,C\,$ is a Cantor set contained in $\,K\,$ then 
$\,\xi+C\,$ is a Cantor set contained in 
$\,{\Bbb R}\setminus K\,$, {\it q.e.d.}
\medskip
From Lemma 1 we immediately derive the following useful observation.
\medskip
(2.2)\quad {\it If $\,B\,$ is a Bernstein set and $\,K\,$ is a proper
subfield of $\,{\Bbb R}\,$ and $\,B\subset K\,$ then $\,K\,$ is a
Bernstein field.}
\medskip
{\it Remark.} In view of (2.2) Bernstein fields can be 
defined in complete analogy with Cantor fields.
While a proper subfield $\,K\,$ of $\,{\Bbb R}\,$ is a {\it Cantor field} 
if and only if {\it $\,K\,$ contains a Cantor set}, 
$\,K\,$ is a {\it Bernstein field} 
if and only if {\it $\,K\,$ contains a Bernstein set}.
\medskip
As usual, for $\,Y\subset{\Bbb R}\,$ let $\,{\Bbb Q}(Y)\,$ 
denote the smallest subfield of $\,{\Bbb R}\,$ containing 
the set $\,Y\,$. 
To cut short, 
a set $\,S\subset{\Bbb R}\,$ is {\it algebraically independent}
if and only if 
for arbitrary $\,n\in{\Bbb N}\,$ 
one cannot find distinct numbers $\,t_1,...,t_n\,$ in $\,S\,$ 
such that $\;p(t_1,...,t_n)=0\;$ holds 
for some nonconstant 
polynomial $\,p(X_1,...,X_n)\;$ in $\,n\,$ indeterminates
with rational coefficients.
It is well-known (see [6]) that this definition 
is equivalent to the condition that 
$\;t\not\in\overline{{\Bbb Q}(S\setminus\{t\})}\;$
for every $\,t\in S\,$.
Since $\,{\Bbb Q}(S)\not={\Bbb R}\,$ for every algebraically independent $\,S\subset{\Bbb R}\,$,
from Lemma~1 and (2.2) we derive
\medskip
(2.3)\quad {\it If $\,B\,$ is a Bernstein set and $\,S\subset{\Bbb R}\,$ is 
algebraically independent and $\,B\subset S\,$ then $\,S\,$ is a Bernstein 
set and $\,{\Bbb Q}(S)\,$ is a Bernstein field.}
\medskip
{\it Remark.} The fact that 
$\,{\Bbb Q}(S)\not={\Bbb R}\,$ for every algebraically independent $\,S\subset{\Bbb R}\,$
is crucial for (2.3) and for what follows. This fact can be verified by 
several arguments. For example it is a nice exercise to verify that 
$\,\sqrt 2\not \in {\Bbb Q}(S)\,$. Moreover,  
$\,{\Bbb Q}(S)\,$ cannot contain any irrational algebraic number.
This can be verified in the following way. 
Firstly, if $\,\theta\,$ is transcendental over a field $\,K\,$
then every element of $\,K(\theta)\setminus K\,$ is transcendental over 
$\,K\,$. (This important basic fact 
is proved in detail in [6] \S 73.) Secondly and consequently
by induction, for every {\it finite} algebraically independent set 
$\,S\subset{\Bbb R}\,$ the set $\,{\Bbb Q}(S)\setminus{\Bbb Q}\,$ cannot contain 
algebraic numbers. Thirdly, this is also true for {\it infinite} $\,S\,$
because for arbitrary $\,A\subset{\Bbb R}\,$ the field 
$\,{\Bbb Q}(A)\,$ is the union of all fields $\,{\Bbb Q}(E)\,$ with $\,E\,$
running through the finite subsets of $\,A\,$. 
\medskip
To cut short, an 
algebraically independent set $\,S\subset{\Bbb R}\,$
is a {\it transcendence base} if and only if for each 
real number $\,x\not\in S\,$ the set 
$\,S\cup\{x\}\,$ is not algebraically independent. 
As an immediate (and important, well-known) consequence 
of this maximality condition, 
$\,{\Bbb R}\,$ must be an algebraic extension
of $\,{\Bbb Q}(S)\,$. Therefore, we easily obtain the following observation.
\smallskip
(2.4) {\it $\,\overline{K}={\Bbb C}\,$ 
for a subfield $\,K\,$ of $\,{\Bbb R}\,$ 
if and only if $\,K\,$ contains 
some transcendence base.}
\medskip                                   
Clearly, $\,|\overline{{\Bbb Q}(T)}|<{\bf c}\,$ 
whenever $\,T\subset{\Bbb R}\,$ and $\,|T|<{\bf c}\,$.
Thus in view of (2.4) and $\,|{\Bbb C}|={\bf c}\,$ 
the size of a transcendence base must be $\,{\bf c}\,$. 
Furthermore, it is well-known 
(and easy to verify) that the fields $\,{\Bbb Q}(S_1)\,$ and $\,{\Bbb Q}(S_2)\,$
are always isomorphic for algebraically independent sets $\,S_1,S_2\,$
of equal size.
Therefore, from (2.3) and (2.4) we derive the following statement 
which gives a hint how to find an appropriate family $\,{\cal B}_1\,$
in Theorem 1.
\medskip
(2.5)\quad {\it If $\,B\,$ is a Bernstein set and $\,T\,$ is 
a transcendence base and $\,B\subset T\,$ then $\,T\,$ is a Bernstein set 
and $\,{\Bbb Q}(T)\,$ is a
Bernstein field and $\,\overline{{\Bbb Q}(T)}={\Bbb C}\,$.} 
\medskip
Thus we obtain appropriate families $\,{\cal C}_1\,$
and $\,{\cal B}_1\,$ in Theorem 1 
by applying the following theorem which is interesting in its own right.
\medskip
{\bf Theorem 2.} {\it There exist two 
families $\;{\cal Y}_1,{\cal Y}_2\;$ 
of transcendence bases with $\;|{\cal Y}_1|=|{\cal Y}_2|=2^{\bf c}\;$ 
such that each member of 
$\,{\cal Y}_1\,$ is a nowhere dense Lebesgue null set containing 
a Cantor set, while each member of $\,{\cal Y}_2\,$
is a Bernstein set, and for $\,i\in\{1,2\}\,$ and for 
distinct $\;R,S\,\in\,{\cal Y}_i\;$ we have
$\,{\Bbb Q}(R)\not\subset{\Bbb Q}(S)\,$ and in particular
$\,{\Bbb Q}(R)\not={\Bbb Q}(S)\,$.}
\medskip
For the proof of Theorem 2 we need the following lemma, which can easily 
be verified, and two propositions we will prove in 
Sections 5 and 6.
\medskip
{\bf Lemma 2.} {\it If $\,B\,$ is a transcendence base and 
$\,\theta\,$ is an algebraic irrational number 
then $\;(B\setminus X)\cup\{\,x\theta\;|\;x\in X\,\}\;$
is a transcendence base for every $\,X\subset B\,$.}
\medskip
{\bf Proposition 2.} {\it There exists an algebraically independent
Cantor set lying in $\,{\Bbb D}\,$.} 
\smallskip
{\bf Proposition 3.} {\it There exist disjoint Bernstein sets 
$\,A,B\,$ such that $\,A\cup B\,$ is an 

algebraically independent Bernstein set.}

\medskip
{\it Remark.} There is a natural analogy between algebraically 
independent subsets of $\,{\Bbb R}\,$ and              
{\it linearly independent} subsets of $\,{\Bbb R}\,$ (as defined in the usual 
way, see [1], [2]), between transcendence bases and {\it Hamel bases}.
A Hamel basis is a basis of the vector space $\,{\Bbb R}\,$ over the field $\,{\Bbb Q}\,$
and hence a {\it maximal} linearly independent set, 
while a transcendence basis is a {\it maximal} algebraically 
independent set. Consequently, {\it both} the existence of Hamel bases 
and the existence of transcendence bases are guaranteed by {\it one}
routine argument applying  Zorn's Lemma.
Notice that  algebraically independent sets must be
linearly independent, but the converse is not true in general.
(For example, it is evident that $\,\{\pi,\pi^2\}\,$ 
is linearly independent but not algebraically independent.
Moreover, if $\,H\,$ is a Hamel basis then $\,H\,$
cannot be algebraically independent because, trivially, $\,{\Bbb Q}(H)={\Bbb R}\,$.)
\bigskip
{\bff 3. Proof of Theorem 2}
\medskip
Fix four subsets $\,S_1,T_1,S_2,T_2\,$ of $\,{\Bbb R}\,$ of size $\,{\bf c}\,$
such that $\;S_1\cap T_1=S_2\cap T_2=\emptyset\;$ 
and both sets $\,S_1\cup T_1\,$ and $\,S_2\cup T_2\,$ 
are transcendence bases and $\,S_2\,$ 
is a Bernstein set while $\,S_1\,$ is a Cantor set 
and $\,S_1\cup T_1\,$ is a subset of $\,{\Bbb D}\,$.
In view of Propositions 2 and 3 
a choice of such sets is possible since each algebraically independent
set is contained in some transcendence base and since 
each Cantor set can be split into two Cantor sets 
and since $\,{\Bbb Q}({\Bbb D})={\Bbb R}\,$.
(We have $\,{\Bbb Q}({\Bbb D})={\Bbb R}\,$ due to 
$\,{\Bbb D}+{\Bbb D}=[0,2]\,$, and for every Cantor set $\,C\,$ both 
$\;C\cap\,]{-\infty,t}]\;$ and $\;C\cap [t,\infty]\;$
are Cantor sets for some real $\,t\not\in C\,$.)
We define for every set $\,U\subset T_i\,$  
\medskip
\centerline{$g_i(U)\;:=\;
S_i\,\cup\,U\,\cup\,\{\,x{\root 3\of 2}\;|\;x\,\in\,T_i\setminus U\,\}\;.$}
\medskip
By virtue of Lemma 2, $\,g_i(U)\,$ is always a transcendence base. 
(In view of (4.1) in the next section we use 
the factor $\,\root 3\of 2\,$ instead of the simpler factor $\,\sqrt 2\,$.)
Trivially $\,g_1(U)\,$ contains a Cantor set.
As a union of three nowhere dense null sets 
$\,g_1(U)\,$ is a nowhere dense null set. 
By (2.3), $\,g_2(U)\,$ is a Bernstein set.
Therefore Theorem 2 is settled by defining 
\smallskip
\centerline{$\;{\cal Y}_i\,:=\,
\{\,g_i(U)\;|\;U\subset T_i\,\}\;\;\;(i\in\{1,2\})\;$}
\smallskip
because in both cases $\,i\in\{1,2\}\,$ the following statement 
is true. 
\medskip
(3.1)\quad {\it For $\;V,W\subset T_i\;$ the equality 
$\;V=W\;$ follows from $\;{\Bbb Q}(g_i(V))\subset{\Bbb Q}(g_i(W))\,$.}
\medskip
In order to verify (3.1) we show a bit more in view of 
the next section and Section 7. 
If $\,Y\subset{\Bbb R}\,$ is algebraically independent, 
then every algebraic number in the field $\,{\Bbb Q}(Y)\,$ is rational
and hence $\,{\root 3\of 2}\not\in {\Bbb Q}(Y)\,$.
Therefore (3.1) is an immediate consequence of the following statement.
\medskip
(3.2)\quad {\it Let  $\,i,j\in\{1,2\}\,$ and 
$\,V\subset T_i\,$ and $\,W\subset T_j\,$.
Furthermore let $\,L\,$ be a subfield of $\,{\Bbb R}\,$
with $\;\root 3\of 2\not\in L\;$
and $\,g_j(W)\subset L\,$. 
Then $\;g_i(V)\subset L\;$ implies  $\;V=W\,$.}
\medskip
In order to verify (3.2) we firstly point out 
that the inclusion $\,g_i(V)\subset L\,$ enforces 
the identity $\,i=j\,$.
Indeed, due to $\,g_i(V)\subset L\,$ and $\,g_j(W)\subset L\,$,
from $\,i\not=j\,$ we derive that 
$\,L\,$ contains both a Cantor set and a Bernstein set 
which is impossible in view of Lemma~1 since  $\,L\not={\Bbb R}\,$.
Thus we may assume $\,i=j\in\{1,2\}\,$ and 
$\,g_i(V)\subset L\,$ and $\,g_i(W)\subset L\,$
and $\,\root 3\of 2\not\in L\,$. 
To conclude the proof of (3.2) by verifying 
$\,V=W\,$ assume indirectly that there is a real 
$\;x\,\in\,(V\setminus W)\cup(W\setminus V)\,$. 
Then the pair $\,(x,x {\root 3\of 2})\,$ lies in
$\,\,g_i(V)\times g_i(W)\,$ or in $\,\,g_i(W)\times g_i(V )\,$.
Consequently the field $\,L\,$ contains
both reals $\,x,x {\root 3\of 2}\,$ and hence also 
the quotient $\,\root 3\of 2\,$ which is a contradiction.
\bigskip
{\bff 4. Proof of Theorem 1}
\medskip
For every subfield $\,K\,$ of $\,{\Bbb R}\,$
let $\,K^*\,$ denote the intersection of all fields 
$\,L\,$ satisfying $\,K\subset L\subset {\Bbb R}\,$ and the property 
that $\,\sqrt{|x|}\in L\,$ for every $\,x\in L\,$.
Of course, $\;K\subset K^*\subset {\Bbb R}\cap\overline{K}\;$
and $\,\sqrt{|x|}\in K^*\,$ for all $\,x\in K^*\,$ or, equivalently, 
$\,K^{**}=K^*\,$.
Alternatively, $\,K^*\,$ is obtained from $\,K\,$ 
by successively adjoining square roots.
Define inductively $\,W_0\,=\,K\,$ 
and $\;W_{n}\,=\,{\Bbb Q}(\{\,\sqrt{|x|}\;|\;x\in W_{n-1}\,\})\;$
for $\,n\in{\Bbb N}\,$ in order to finally obtain 
$\;K^*=\bigcup_{n=1}^\infty W_n\,$. (As usual in number theory, 
$\,0\not\in{\Bbb N}\,$.)
\medskip
A fortiori 
$\,\overline{{\Bbb Q}(T)^*}={\Bbb C}\,$ 
for every transcendence base $\,T\,$, 
where $\,{\Bbb Q}(T)^*\not={\Bbb R}\,$ due to the following observation.
\smallskip
(4.1)\quad {\it If $\,Y\subset {\Bbb R}\,$ is an algebraically independent set
then $\,\root 3\of 2\not\in{\Bbb Q}(Y)^*\,$.}
\medskip
Obviously (4.1) is a consequence of the following two statements. 
\smallskip
(4.2)\quad {\it If $\,K\,$ is a subfield of $\,{\Bbb R}\,$ then
every number in $\,K^*\setminus K\,$ 
is algebraic over $\,K\,$ of degree $\,2^n\,$ with $\,n\in{\Bbb N}\,$.}
\smallskip
(4.3)\quad {\it For every algebraically independent set $\,Y\subset{\Bbb R}\,$
the degree of $\,\root 3\of 2\,$ 
over the field $\,{\Bbb Q}(Y)\,$ equals $\,3\,$.}
\medskip
We obtain (4.2) from the well-known fact that 
$\,K^*\,$ is the field of all reals which are 
{\it constructible} over $\,K\,$ by a finite sequence of 
{\it ruler and compass} constructions.
\smallskip
In order to verify (4.3) let 
$\,Y\subset{\Bbb R}\,$ be algebraically 
independent and consider distinct $\;y_1,...,y_n\in Y\;$
for arbitrary $\,n\in{\Bbb N}\,$.
Naturally the field $\,{\Bbb Q}(y_1,...,y_n)\,$ 
is isomorphic with the quotient field $\,{\Bbb Q}(X_1,...,X_n)\,$
of the polynomial ring $\,{\Bbb Z}[X_1,...,X_n]\,$ which is 
a unique factorization domain and in which $\,2\,$ is irreducible.
Consequently, by Eisenstein's Criterion the polynomial $\,X^3-2\,$ 
is irreducibel over the field $\,{\Bbb Q}(X_1,...,X_n)\,$. 
Thus $\,X^3-2\,$ is irreducibel over $\,{\Bbb Q}(y_1,...,y_n)\,$ and hence
$\,X^3-2\,$ is the minimal polynomial of 
$\,\root 3\of 2\,$ over the field $\,{\Bbb Q}(Y)\,$ and (4.3) is proved.
\medskip\smallskip
Now we take the transcendence bases $\,g_i(\cdot)\,$ from 
the previous section and define
\medskip
\centerline{$\;{\cal C}_2\,:=\,\{\,{\Bbb Q}(g_1(U))^*\;|\;U\subset T_1\,\}\;$
\quad and \quad
$\;{\cal B}_2\,:=\,\{\,{\Bbb Q}(g_2(U))^*\;|\;U\subset T_2\,\}\,$.}
\medskip
Then all fields in $\,{\cal C}_2\,$ resp.~$\,{\cal B}_2\,$
are Cantor fields resp.~Bernstein fields. 
Furthermore, concluding the proof of Theorem 1, the following statement shows 
that $\,|{\cal C}_2|=|{\cal B}_2|=2^{\bf c}\,$ and 
the fields in $\,{\cal C}_2\cup{\cal B}_2\,$
are incomparable.
\medskip
(4.4)\quad {\it For $\,i,j\in\{1,2\}\,$ and 
$\,V\subset T_i\,$ and $\,W\subset T_j\,$
the field $\,{\Bbb Q}(g_i(V))^*\,$ can be embedded into the field
$\,{\Bbb Q}(g_j(W))^*\,$ only if $\,V=W\,$.}
\medskip
The following lemma shows that if 
$\,{\Bbb Q}(g_i(V))^*\,$ can be embedded into 
$\,{\Bbb Q}(g_j(W))^*\,$ then already 
$\,{\Bbb Q}(g_i(V))^*\subset {\Bbb Q}(g_j(W))^*\,$ holds.
Consequently and in view of (4.1) we obtain (4.4) from (3.2) 
with $\,L:={\Bbb Q}(g_j(W))^*\,$.
\medskip\smallskip
{\bf Lemma 3.} {\it If $\,K\,$ is a subfield of $\,{\Bbb R}\,$
with $\,K^*=K\,$ and $\,\varphi\,$ is a monomorphism
from $\,K\,$ to $\,{\Bbb R}\,$ then $\,\varphi(x)=x\,$ for 
all $\,x\in{\Bbb R}\,$.}
\medskip
{\it Proof.} It goes without saying that
$\,\varphi(r)=r\,$ for all $\,r\in{\Bbb Q}\,$.
Due to $\,K^*=K\,$ we have 
$\;\sqrt{|a|}\in K\;$ for every $\,a\in K\,$.
Now let $\,x\in K\,$ and $\,r\in{\Bbb Q}\,$.
If $\,r<x\,$ then $\,\varphi(x)-\varphi(r)=\varphi(x-r)=
\varphi(\sqrt{x-r})^2>0\,$, and if 
$\,r>x\,$ then $\,\varphi(r)-\varphi(x)=\varphi(r-x)=
\varphi(\sqrt{r-x})^2>0\,$. Thus $\,\varphi\,$ is strictly increasing 
and, since 
$\,\varphi\,$ restricted to the domain  $\,{\Bbb Q}\,$ is the identity,
for all $\,r,s\in{\Bbb Q}\,$ and $\,x\in K\,$ 
the implication $\;r<x<s\;\Longrightarrow\;r<\varphi(x)<s\;$ holds
and this enforces $\;\varphi(x)=x\;$ for all  $\,x\in K\,$,  {\it q.e.d.}
\bigskip
{\bff 5. Proof of Proposition 2}
\medskip
As usual, ZF means ZFC set theory minus the Axiom of Choice.
In this section we are going to prove Proposition 2 without 
applying the Axiom of Choice. This has the benefit that the 
following trivial consequence of Theorem 1 turns out to be a theorem in ZF.
\medskip
{\bf Corollary 1.} {\it There exist two 
families $\,{\cal C}_3,{\cal C}_4\,$
equipollent with the power set of $\,{\Bbb R}\,$ 
such that $\,{\cal C}_3\,$ consists of isomorphic Cantor fields 
while $\,{\cal C}_4\,$ consists of incomparable Cantor fields.} 
\medskip
{\it Remark.} There are two reasons why Theorem 1 is not a theorem in ZF.
Firstly, the existence of a proper subfield $\,K\,$ of $\,{\Bbb R}\,$ 
satisfying $\,\overline{K}={\Bbb C}\,$ is unprovable in ZF
because of (2.4) and the well-known fact that the existence 
of transcendence bases is unprovable in ZF. Secondly, 
also the existence of a Bernstein set (and in particular 
the existence of a Bernstein field) cannot be derived from the axioms 
of ZF only. Additionally, equations like $\,|{\cal F}|=2^{\bf c}\,$
need to be interpreted in a specific way in ZF since 
$\,|X|\,$ is not defined in ZF for arbitrary sets $\,X\,$.
(But notice that sets like $\,{\Bbb R},{\Bbb C}\,$ and $\,{\Bbb D}\,$
are well-defined objects in ZF. Notice also that in ZF the field
$\,\overline{K}\cap{\Bbb R}\,$ is well-defined for each subfield  $\,K\,$ of $\,{\Bbb R}\,$
constructed in ZF.)

\medskip
If $\,C\,$ is an
algebraically independent Cantor set in ZF then 
a ZF-proof of Corallary 1 can easily be obtained by adopting 
parts of the ZFC-proof of Theorem 1. 
Simply split $\,C\,$ constructively into two disjoint 
Cantor sets $\,C_1,C_2\,$ and define
the sets $\,S_1\,$ and $\,T_1\,$ in the proof of 
Theorem 1 via $\,S_1=C_1\,$ and $\,T_1=C_2\,$. 
Then, never considering transcendence bases,
$\;C_1\,\cup\,U\,\cup\,\{\,x{\root 3\of 2}\;|\;x\,\in\,C_2\setminus U\,\}\;$
is algebraically independent
for all $\,U\subset C_2\,$ and this is enough.
Alternatively, there is a direct and much shorter ZF-proof of Corollary~1
as follows. With $\,C_1,C_2\,$ as above put 
$\;{\cal C}_3\,:=\,\{\,{\Bbb Q}(C_1\cup X)\;|\;X\subset C_2\,\}\;$ and 
$\;{\cal C}_4\,:=\,\{\,{\Bbb R}\cap\overline{{\Bbb Q}(C_1\cup X)}\;|\;X\in {\cal F}\}\;$
where $\,{\cal F}\,$ is a family equipollent with the power set of $\,{\Bbb R}\,$ 
such that $\,X\subset C_2\,$ for every $\,X\in{\cal F}\,$ 
and $\,X\not\subset Y\,$ whenever $\,X,Y\in{\cal F}\,$ are distinct.
(Such a family $\,{\cal F}\,$ can easily be constructed in ZF. 
For example take a bijection $\,f\,$ from $\,{\Bbb R}\,$
onto $\,C_2\,$ and put 
$\;{\cal F}\,:=
\,\{\,f(T\cup ([2,3]\setminus(2\!+\!T)))\;|\;T\subset [0,1]\,\}\,$.)
The family $\,{\cal C}_3\,$  obviously fits
and the family $\,{\cal C}_4\,$ does the job by virtue of Lemma 3.
\medskip
There is a more algebraic way and a more topological way to prove
Proposition 2 in ZF. The algebraic way is far from being elementary  
since an appropriate modification 
of the famous {\it von Neumann numbers} (see [5]) is used.
\smallskip
As usual, $\,[\xi]\,$ denotes the largest integer
$\,k\leq\xi\,$ for $\,\xi\in{\Bbb R}\,$. We put 
$\;\psi(x,n):=2^{n^2}-2^{[nx]}\;$ 
and observe that 
$\;0<\psi(x,n)<\psi(x,m)\;$ whenever $\,n,m\in{\Bbb N}\,$ and
$\,0<x<n<m\,$. Consequently, 
\smallskip
\centerline{$\;\sigma(x)\;:=\;2\cdot\!\!\!\sum\limits_{n>x} 3^{-\psi(x,n)}\;$}
\smallskip
defines a function from $\;]0,\infty[\;$ into $\,{\Bbb D}\,$ for alle $\,x\in{\Bbb R}\,$.
(Notice that the von Neumann numbers are not contained in $\,{\Bbb D}\,$.)
In the same way as carried out in [5] one can verify that 
$\;\sigma(x)\not\in\overline{{\Bbb Q}(\sigma(]{0,x}[))}\;$
for every $\,x>0\,$. Hence $\,\sigma\,$ is 
injective and the set $\,\sigma(]0,\infty[)\,$ is algebraically independent. 
It goes without saying that $\,\sigma\,$ is continuous at each positive
irrational. Now consider the uncountable partition $\,{\cal P}\,$ of $\,{\Bbb D}\,$
defined in the proof of (2.1) and obviously well-defined in ZF.
Then $\;{\cal P}^*\,:=\,\{\,X\in{\cal P}\;|\;X\cap{\Bbb Q}\not=\emptyset\,\}\;$
is a countable subfamily of $\,{\cal P}\,$ and hence the family 
$\,{\cal P}\setminus{\cal P}^*\,$ is not empty. Hence we can select 
an element $\,C\,$ from this family. 
Automatically, $\,C\,$ is a Cantor set
containing only irrational numbers.
Consequently, $\,\sigma(C)\,$ is a Cantor set lying in an 
algebraically independent set and this concludes 
the ZF-proof of Proposition 2.
\medskip\smallskip
Alternatively we now present an elementary ZF-proof of Proposition 2.
(Notice that $\,{\Bbb D}\,$ is well-defined in ZF).
In doing so we show a bit more, namely that 
the following statement (which cannot be verified using 
some modification of the von
Neumann numbers) is a theorem in ZF.
\medskip
(5.1)\quad{\it Every Cantor set $\,D\,$ contains an algebraically independent 
Cantor set.}
\medskip
First of all recall the standard construction of 
the Cantor ternary set 
\smallskip
\centerline{${\Bbb D}\;=\;\bigcap\limits_{k=1}^\infty
\bigcup\limits_{j=1}^{2^k}\;{\cal J}_{k,j}$}
\smallskip
where $\,{\cal J}_{1,1}=[0,{1\over 3}]\,$ and 
$\,{\cal J}_{1,2}=[{2\over 3},1]\,$ and 
$\,{\cal J}_{2,1}=[0,{1\over 9}]\,$ and 
$\,{\cal J}_{2,2}=[{2\over 9},{1\over 3}]\,$ and
$\,{\cal J}_{2,3}=[{2\over 3},{7\over 9}]\,$ and
$\,{\cal J}_{2,4}=[{8\over 9},1]\,$ and so on.
\medskip
For $\,n\in{\Bbb N}\,$ consider the ring 
$\;R_n={\Bbb Z}[X_1,...,X_n]\;$
of all integral polynomials in $\,n\,$ indeterminates.
Of course, $\,R_n\,$ is a countable set.
For every $\,n\in{\Bbb N}\,$ let $\,g_n\,$ be a bijection 
from $\,{\Bbb N}\,$ onto $\,R_n\setminus{\Bbb Z}\,$.
For $\,(n,m)\in{\Bbb N}^2\,$ let $\,A(n,m)\,$ denote 
the set of all $n$-tuples $\,(t_1,...,t_n)\,$ of reals
that are annihilated by the polynomial  $\,g_n(m)\,$.
Naturally, $\,A(n,m)\,$ is a {\it closed} subset of the space $\,{\Bbb R}^n\,$.
Furthermore, $\,A(n,m)\,$ is {\it nowhere dense} because
for an integral polynomial 
$\,p(X_1,...,X_n)\,$ the equation 
$\,p(t_1,...,t_n)=0\,$ holds for all points $\,(t_1,...,t_n)\,$ 
in a nonempty open subset of $\,{\Bbb R}^n\,$ only if 
$\,p(X_1,...,X_n)\,$ is the zero polynomial.
\medskip
For $\,a\leq b\,$ put $\,\lambda([a,b]):=b-a\,$.
For every $\,k\in{\Bbb N}\,$ choose $\,2^k\,$ pairwise disjoint
compact intervals
$\;{\cal I}_{k,j}\;(1\leq j\leq 2^k)\;$
such that $\;{\cal I}_{m,j}\;\supset\;
{\cal I}_{m+1,2j}\cup{\cal I}_{m+1,2j-1}\;$ 
whenever $\,m,j\in{\Bbb N}\,$ and $\,j\leq 2^m\,$ and that 
$\;\lim_{k\to\infty}\max\,\{\,\lambda({\cal I}_{k,j})\;|\;1\leq j\leq 2^k\,\}
\,=0\,$.
Certainly, in each parallelepiped of positive volume
there lies a parallelepiped of positive volume
disjoint with the closed, and nowhere dense point set 
$\,A(n,m)\,$. Therefore, step by step it can be accomplished 
that the following condition $\,{\bf B}[m]\,$ holds for all $\,m\in{\Bbb N}\,$ 
as well.
\smallskip
$({\bf B}[m])\;$ {\it $\,A(n,m)\cap\prod_{i=1}^n\varphi(i)\,=\,\emptyset\;$
whenever $\,n\in{\Bbb N}\,$ and $\,\varphi\,$ is an
injection from $\;\{1,...,n\}\;$
into $\;\{{\cal I}_{m,1},...,{\cal I}_{m,2^m}\}\,$.}
\medskip
In comparison with the standard construction of the ternary Cantor set 
it is clear that 
\smallskip
\centerline{$C\;:=\;\bigcap\limits_{k=1}^\infty
\bigcup\limits_{j=1}^{2^k}\;{\cal I}_{k,j}$}
\smallskip
is a Cantor set. (Alternatively it is plain that 
$\,C\,$ is compact, dense in itself and does not contain isolated points.)
We claim that $\,C\,$ is algebraically independent.
\medskip
Assume indirectly that $\,C\,$ is not algebraically independent.
Then we can find distinct numbers $\;t_1,t_2,...,t_n\;$ in $\,C\,$
such that the $n$-tuple $\,(t_1,t_2,...,t_n)\,$ is annihilated by some 
non-zero polynomial $\,p\,$ in the ring $\,{\Bbb Z}[X_1,...,X_n]\,$.
Trivially, this $n$-tuple is also annihilated 
by the polynomial $\,k\cdot p(X_1,...,X_n)\,$ 
for every $\,k\in{\Bbb N}\,$. 
Hence there is an infinite set $\,M\subset {\Bbb N}\,$
such that for all $\,m\in M\,$
the point $\,(t_1,...,t_n)\,$ lies in the set $\,A(n,m)\,$.
Now choose $\,m\in M\,$ sufficiently large such that 
each one of the $\,2^m\,$ intervals $\,{\cal I}_{m,j}\,$
contains at most one of the numbers $\,t_1,...,t_n\,$.
This provides us with an injection $\,\varphi\,$ 
from $\;\{1,...,n\}\;$ into  $\;\{{\cal I}_{m,1},...{\cal I}_{m,2^m}\}\;$
such that the point 
$\,(t_1,...,t_n)\,$ 
lies in the parallelepiped
$\;\prod_{i=1}^n\varphi(i)\,$. 
This is a contradiction to $\,{\bf B}[m]\,$
since $\,(t_1,...,t_n)\in A(n,m)\,$.
\medskip
So we have proved that $\,C\,$ is algebraically independent.
In doing so also Proposition 2 is settled because  
in choosing the intervals $\,{\cal I}_{k,j}\,$
we can take care that always 
$\,{\cal I}_{k,j}\subset{\cal J}_{k,j}\,$, which trivially implies 
$\,C\subset {\Bbb D}\,$.
Moreover, an appropriate choice of the intervals $\,{\cal I}_{k,j}\,$
is also possible to accomplish 
$\,C\subset D\,$ for an arbitrary Cantor set
$\,D\,$ and hence (5.1) is settled. This is indeed true 
in view of the following proposition 
about Cantor sets and because 
$\;f\big(\bigcap_k\bigcup_j\;{\cal I}_{k,j}\big)\,=\,
\bigcap_k\bigcup_j\;f({\cal I}_{k,j})\;$
and $\;f([u,v])=[f(u),f(v)]\;$ 
for every increasing bijection $\;f:\,{\Bbb R}\to{\Bbb R}\,$.
\medskip
{\bf Proposition 4.} {\it A nonempty set $\,C\subset{\Bbb R}\,$ 
is compact, totally disconnected and dense in itself
if and only if $\,C=f({\Bbb D})\,$ 
for some increasing bijection $\,f\,$ from $\,{\Bbb R}\,$ onto $\,{\Bbb R}\,$.}
\medskip
{\it Proof.} Naturally, every   
strictly increasing function $\,f\,$ from $\,{\Bbb R}\,$ onto  $\,{\Bbb R}\,$
is a homeomorphism. Hence as the topological space $\,{\Bbb D}\,$
the set $\,f({\Bbb D})\,$ is compact, totally disconnected and dense in itself.
Conversely let $\,D_1=C\,$ similarly as $\,D_0={\Bbb D}\,$ be a 
compact, totally disconnected and dense-in-itself 
subspace of $\,{\Bbb R}\,$. In the following, $\,i\in\{0,1\}\,$.
\smallskip
The set $\,D_i\,$ has a maximum and a minimum due to compactness.
Let $\,a_i=\min D_i\,$ and $\,b_i=\max D_i\,$.
(In particular, $\,a_0=0\,$ and $\,b_0=1\,$.)
Since the set $\;[a_i,b_i]\setminus D_i\;$
is open, we can define 
a family $\,{\cal U}_i\,$ of pairwise
disjoint open intervals such that
$\;\bigcup{\cal U}_i\subset[a_i,b_i]\;$ 
and $\;D_i\,=\,[a_i,b_i]\setminus\bigcup{\cal U}_i\,$.
\smallskip
In a natural way the family $\,{\cal U}_i\,$ is linearly ordered 
by declaring 
$\,U\in{\cal U}_i\,$ {\it smaller than} $\,V\in{\cal U}_i\,$ 
if and only if $\;\sup U\,\leq\,\inf V\,$. 
Since $\,D_i\,$ is totally disconnected, 
$\,\bigcup{\cal U}_i\,$ is dense in $\,[a_i,b_i]\,$.
And since $\,D_i\,$ has no isolated point, 
between two intervals from $\,{\cal U}_i\,$ there always lie 
infinitely many intervals from $\,{\cal U}_i\,$.
Finally, the linearly orderd set 
$\,{\cal U}_i\,$ has neither a smallest nor a largest element.
\smallskip
Taking all these observations into account, 
by virtue of a well-known classical theorem by Cantor, the linearly ordered
set  $\,{\cal U}_i\,$ is order-isomorphic to 
the naturally ordered set $\,{\Bbb Q}\,$.
In particular the two families 
$\,{\cal U}_0\,$ and $\,{\cal U}_1\,$ are order-isomorphic.
Since for $\,u<v\,$ and $\,r<s\,$ the interval $\;]u,v[\;$
can easily be mapped onto 
the interval $\;]r,s[\;$ by a strictly increasing function,
with the help of an order-isomorphism from 
$\,{\cal U}_0\,$ onto $\,{\cal U}_1\,$
we obtain an increasing bijection $\,g\,$
from $\,\bigcup{\cal U}_0\,$ onto $\,\bigcup{\cal U}_1\,$.
This function $\,g\,$ can be extented to a continuous function 
$\,f\,$ from $\,{\Bbb R}\,$ to  $\,{\Bbb R}\,$ by defining 
$\;f(t)\,:\,=a_1+t\;$ for $\,t\leq 0\,$ 
and $\;f(t)\,:=\,b_1+t-1\;$ for $\,t\geq 1\,$ and 
\smallskip
\centerline{$\;f(t)\,:=\,
\sup\,\{\,g(x)\;|\;x\in\bigcup{\cal U}_0\;\land\;x\leq t\,\}\;$}
\smallskip
for $\,0<t<1\,$. Obviously $\,f\,$ is strictly increasing and surjective 
and $\,f({\Bbb D})=C\,$, {\it q.e.d.}
\medskip\smallskip
{\it Remark.} Proposition 4 can be regarded as a 
{\it purely algebraic} characterization of Cantor sets.
Because firstly in the definition of $\,{\Bbb D}\,$ 
the series $\,\sum a_n3^{-n}\,$ need not refer to the 
topological concept {\it convergency} but can be regarded 
as the depiction of a sequence of digits.
And secondly, {\it monotonicity} 
of real functions can be characterized 
purely algebraically because a mapping $\,f\,$ 
from $\,{\Bbb R}\,$ into $\,{\Bbb R}\,$ is increasing if and only if 
for arbitrary $\,u,v\in{\Bbb R}\,$
the equation $\;(u-v)(f(u)-f(v))=x^2\;$ has a real solution $\,x\,$.
\medskip
{\it Remark.} The existence of an algebraically independent Cantor set, 
which is a trivial consequence of (5.1), can be obtained by applying 
the Kuratowski-Mycielski theorem [2] 19.1. However, this theorem is 
not elementary (since arguments using Vietoris topologies occur) and it 
does not imply (5.1), even not in ZFC.
On the other hand, our proof of (5.1) is elementary 
and, since it is carried out in ZF, constructive as well.
In view of [2] 19.2(i) there exists a Hamel basis which contains 
a Cantor set. Even more, by applying (5.1) {\it there exist a Cantor set $\,C\,$
and a transcendence basis $\,T\,$ and a Hamel basis $\,H\,$ such that
$\;C\subset T\subset H\subset {\Bbb D}\,$.} (Inevitably, $\,T\not= H\,$ 
since $\,H\,$ cannot be algebraically independent.)
\bigskip
{\bff 6. Proof of Proposition 3}
\medskip
Let $\,{\cal D}\,$ be the family of all Cantor sets.
We have $\,|{\cal D}|={\bf c}\,$ since  
every set in $\,{\cal D}\,$ is compact and $\,{\Bbb R}\,$ 
contains precisely $\,{\bf c}\,$ closed sets and since the translate 
$\;x+{\Bbb D}\;$ is obviously a Cantor set for each one of the $\,{\bf c}\,$ 
real numbers $\,x\,$.
\smallskip
In the following we regard the cardinal 
$\,{\bf c}\,$ as an initial ordinal number and use the ordinal numbers 
$\,\alpha<{\bf c}\,$ for indexing real numbers and sets of real numbers.
(Keep in mind the crucial estimate 
$\;|\{\,\beta\;|\;\beta<\alpha\,\}|<{\bf c}\;$ whenever
$\,\alpha<{\bf c}\,$.) 
Since $\;|\{\,\alpha\;|\;\alpha<{\bf c}\,\}|={\bf c}\,$,
we can write $\;{\cal D}\;=\;\{\,D_\alpha\;|\;\alpha<{\bf c}\,\}\;$ 
(and we do not care whether the mapping $\,\alpha\mapsto D_\alpha\,$
is injective or not).
Let $\,\rho\,$ be a choice function defined 
on the nonempty subsets of $\,{\Bbb R}\,$.
Thus $\,\rho(X)\in X\,$ 
for  every nonempty set $\,X\subset{\Bbb R}\,$. 
\smallskip
Now, by induction we define real numbers $\,x_\alpha\,$ and $\,y_\alpha\,$
for all $\,\alpha<{\bf c}\,$.
Assume that for $\,\xi<{\bf c}\,$ real numbers 
$\,x_\alpha\,$ and $\,y_\alpha\,$ are already defined 
for all $\,\alpha<\xi\,$. (This assumption is vacuous if
$\,\xi=0\,$.)
Then define necessarily distinct 
real numbers $\,x_\xi\,$ and $\,y_\xi\,$ in two steps.
Put 
\medskip  
\centerline{$x_\xi\;:=\;\rho\big(\,D_\xi\,\setminus\,
\overline{{\Bbb Q}(\{\,x_\alpha \;|\;\alpha<\xi\,\}
\cup\{\,y_\alpha \;|\;\alpha<\xi\,\})}\,\big)$}
\medskip
and then, with respect to this definition of $\,x_\xi\,$, put
\medskip  
\centerline{$y_\xi\;:=\;\rho\big(\,D_\xi\,\setminus\,
\overline{{\Bbb Q}(\{\,x_\alpha \;|\;\alpha\leq \xi\,\}
\cup\{\,y_\alpha \;|\;\alpha<\xi\,\})}\,\big)\;.$}
\medskip              
Both definitions are correct since the choice function $\,\rho\,$ 
is always applied to a nonempty set. 
Indeed, for arbitrary $\;Y,T\subset {\Bbb R}\;$ 
the set
$\;Y\setminus\overline{{\Bbb Q}(T)}\;$ certainly is nonempty provided
that $\,|Y|={\bf c}\,$ and $\,|T|<{\bf c}\,$ because 
the field $\;\overline{{\Bbb Q}(T)}\;$ is either countable or
equipollent with $\,T\,$.
So in this way we have 
defined real numbers 
$\,x_\alpha\,$ and $\,y_\alpha\,$ for each $\,\alpha<{\bf c}\,$ 
such that $\,x_\alpha\not=x_\beta\,$ and $\,y_\alpha\not=y_\beta\,$ 
whenever $\,\alpha\not=\beta\,$
and that  $\;\{\,x_\alpha \;|\;\alpha< {\bf c}\,\}\cap
\{\,y_\alpha \;|\;\alpha< {\bf c}\,\}\,=\,\emptyset\,$.
Moreover, the specific choices of the real numbers $\,x_\alpha,y_\alpha\,$
immediately yield 
\medskip
(6.1)\quad {\it The set  
$\;\{\,x_\alpha \;|\;\alpha< {\bf c}\,\}\cup
\{\,y_\alpha \;|\;\alpha< {\bf c}\,\}\;$
is algebraically independent.}
\medskip
Finally, in view of (6.1) and (2.3) 
Proposition 3 is settled with the definitions 
$\;A\,:=\,\{\,x_\alpha \;|\;\alpha< {\bf c}\,\}\;$ and 
$\;B\,:=\,\{\,x_\alpha \;|\;\alpha< {\bf c}\,\}\;$ 
because $\,A\,$ and $\,B\,$ are disjoint and,
due to $\;\{x_\alpha,y_\alpha\}\subset D_\alpha\;$ for every
$\,\alpha<{\bf c}\,$, $\,A\,$ meets every Cantor set
and $\,B\,$ meets every Cantor set. 
\bigskip
{\bff 7. Extremely large fields}
\medskip
Due to the condition $\,\overline{K}={\Bbb C}\,$ the fields depicted 
in Theorem 1 are rather large. But in a certain sense they are not extremely 
large. Let $\,K\,$ be a subfield of $\,{\Bbb R}\,$.
Naturally, $\,{\Bbb R}\,$ is a vector space over the field 
$\,K\,$. As usual, $\;[{\Bbb R}:K]\;$ denotes the dimension of the 
$K$-vector space $\,{\Bbb R}\,$. Naturally,  
in the nontrivial case $\,K\not={\Bbb R}\,$ 
the dimension $\,[{\Bbb R}:K]\,$ is always an infinite cardinal number 
not greater than $\,{\bf c}=|{\Bbb R}|\,$.
\medskip
For every field $\,K\,$ in the family 
$\,{\cal C}_1\cup{\cal C}_2\cup{\cal B}_1\cup{\cal B}_2\,$
defined in the proof of Theorem 1 we have 
$\,[{\Bbb R}:K]={\bf c}\,$ because 
if $\,T\,$ is a transcendence base then 
$\,[{\Bbb R}:{\Bbb Q}(T)]={\bf c}\,$ (see [3] Theorem 4)
and we have 
$\,[{\Bbb R}:{\Bbb Q}(T)^*]={\bf c}\,$ as well. 
(Of course, $\,[{\Bbb R}:{\Bbb Q}(T)^*]={\bf c}\,$ implies $\,[{\Bbb R}:{\Bbb Q}(T)]={\bf c}\,$
and $\,[{\Bbb R}:{\Bbb Q}(T)^*]={\bf c}\,$ can be proved  by 
verifying that the $\,{\bf c}\,$ reals $\;\root 3\of{t}\;(t\in T)\;$
are linearly independent vectors in the vector space $\,{\Bbb R}\,$ 
over the field $\,{\Bbb Q}(T)^*\,$.)
\medskip
One may regard a subfield $\,K\,$ of $\,{\Bbb R}\,$ 
{\it larger than} a subfield $\,L\,$ of $\,{\Bbb R}\,$ 
when $\,[{\Bbb R}:K]<[{\Bbb R}:L]\,$. In this sense the following theorem 
provides extremely many mutually non-isomorphic 
{\it extremely large} Cantor fields and Bernstein fields, 
respectively. Notice that (as explained in [3])
we must have $\,\overline{K}={\Bbb C}\,$ for a subfield $\,K\,$ of $\,{\Bbb R}\,$
satisfying the condition $\,[{\Bbb R}:K]=\aleph_0\,$
\medskip
{\bf Theorem 3.} {\it In Theorem 1 it can be accomplished 
that $\,[{\Bbb R}:K]=\aleph_0\,$ for every field $\,K\,$
in the family $\,{\cal C}_2\cup{\cal B}_2\,$.} 
\medskip
In order to settle Theorem 3 we consider the families 
$\,{\cal C}_2,{\cal B}_2\,$ defined in Section 4, 
expand each field $\,K\,$ in 
$\,{\cal C}_2\cup{\cal B}_2\,$ to an appropriate 
proper subfield $\,\hat K\,$, and replace  
$\,{\cal C}_2\,$ and $\,{\cal B}_2\,$ 
by $\;\{\,\hat K\;|\;K\in{\cal C}_2\,\}\;$ and
$\;\{\,\hat K\;|\;K\in{\cal B}_2\,\}\,$. 
In view of (3.2) this is enough provided that 
for each subfield $\,K\,$ of $\,{\Bbb R}\,$
with $\,\root 3\of 2\not\in K\,$ we can find 
a subfield $\,\hat K\,$ of $\,{\Bbb R}\,$            
with $\,\root 3\of 2\not\in \hat K\,$ and $\,\hat K\supset K\,$
and $\,[{\Bbb R}:\hat K]=\aleph_0\,$.
\medskip
So let $\,K\,$ be a subfield of $\,\R\,$
with $\,\root 3\of 2\not\in K\,$.
Then the family $\,{\cal F}\,$ of all subfields 
$\,L\,$ of $\,\R\,$ with $\,\root 3\of 2\not\in L\,$
and $\;L\supset K\;$ is not empty.
Clearly, $\,\bigcup{\cal G}\in{\cal F}\,$ for 
every chain $\,{\cal G}\,$ of fields in $\,{\cal F}\,$.
Consequently, by applying Zorn's Lemma,
the partially ordered family $\,({\cal F},\subset)\,$ has 
a maximal element $\,\hat K\,$. 
Such a field $\,\hat K\,$ is a subfield $\,\K\,$ of $\,\R\,$ 
satisfying for the irrational number $\,\theta=\root 3\of 2\,$ 
the property depicted in the following proposition.
This concludes the proof of Theorem 3. 
\mp
{\bf Proposition 5.} {\it Let $\,\theta\in\R\,$ be irrational 
and let $\,\K\,$ be a subfield of $\,\R\,$ such that  
$\,\theta\not\in \K\,$
but $\,\theta\in L\,$ for every field 
$\,L\,$ with  $\,\K\subset L\subset\R\,$ and $\,L\not=\K\,$.
Then $\,[\R:\K]=\aleph_0\,$.}
\mp
The best way to verify Proposition 5 is to 
apply methods from infinite-dimensional 
Galois as carried out in [4]. To keep this article self-contained, 
we write down the following proof of Proposition 5 which is 
identical with the proof of [4] Proposition 3.
\mp
Let $\,\theta,\K\,$ be as depicted in Proposition 5.
First it is plain that $\,{\Bbb R}\,$ lies algebraically over $\,\K\,$.
Hence $\,{\Bbb C}\,$ is an algebraic Galois extension
of $\,\K\,$. Let $\,G\,$ be its Galois group.
Let $\,\theta'\not=\theta\,$ be a
$\K$-conjugate of $\,\theta\,$
and choose $\,\tau\in G\,$ such that
$\,\tau(\theta)=\theta'\,$.
Further, let $\,\sigma\,$ be the complex conjugation
$\;z\mapsto \overline z\,$. 
Then $\,\sigma\in G\,$ since $\,\K\subset{\Bbb R}\,$.
Now let $\,H\,$ be the smallest closed subgroup of
$\,G\,$ with $\;\sigma,\tau\in H\,$.
Then $\,\K\,$ must be the fixed field
$\;F\,=\,\{\,z\in{\Bbb C}\;\,|\,\;\forall\,\varphi\in H:\;
\varphi(z)=z\,\}\;$ because $\,\sigma\in H\,$ implies
$\;\K\subset F\subset{\Bbb R}\;$ and $\,F\not=\K\,$
would contradict the maximality of $\,\K\,$
since $\,\theta\not\in F\,$. Consequently, $\,G=H\,$.
Therefore $\,G\,$ is the closure of a finitely
generated subgroup and hence it is easy to verify 
that the family $\,{\cal N}\,$ of all
closed normal subgroups $\,N\,$
of $\,G\,$ with finite index $\,[G:N]\,$
is countable. Since the family 
$\,{\cal L}\,$ of all normal field extensions
$\,L\subset{\Bbb C}\,$ of $\,\K\,$
with finite degree $\,[L:\K]\,$
is equipollent to $\,{\cal N}\,$ because of the Galois correspondence,
$\,{\cal L}\,$ is countable too.
Since every finite extension lies under
a suitable normal finite extension, 
$\;{\Bbb C}\,=\,\bigcup{\cal L}\,$.
Now let $\,B\,$ be a basis of the vector space $\,{\Bbb C}\,$
over $\,\K\,$. Then $\;B\cap L\;$ is finite for 
all $\;L\in{\cal L}\;$ and hence $\;B=B\cap{\Bbb C}\;$
is a countable union of the finite sets
$\;B\cap L\,\,(L\in{\cal L})\,$. Consequently,
$\;\aleph_0\,=\,|B|\,=\,[{\Bbb C}:\K]\,\geq\,[{\Bbb R}:\K]\,\geq\,\aleph_0\,$.

\bigskip\bigskip\bigskip
{\bff References} 
\medskip\smallskip
[1] Ciesielski, K.: {\it Set Theory for the Working Mathematician.} 
Cambridge 1997.
\medskip
[2] Kechris, A., {\it Classical Descriptive Set Theory}, Springer 1995.

\medskip
[3] Kuba, G.: {\it šber den Grad unendlicher K"rpererweiterungen.}
Mitt.~Math.~Seminar

\rightline{Univ.~Hamburg {\bf 26}, 117-124 (2007)} 
\smallskip
[4] Kuba, G.: {\it Counting fields of complex numbers.}
Am.~Math.~Monthly {\bf 116}, No.~6, 

\rightline{541-547 (2009)}
\smallskip
[5] von Neumann, J.: {\it Ein System algebraisch unabh{\"a}ngiger Zahlen.}

\rightline{Math.~Ann.~{\bf 99} (1928) 134-141.}
\smallskip
[6] van der Waerden, B.L.: {\it Algebra I.} Springer-Verlag,  1971.

\bigskip\bigskip\medskip
{\sl Author's address:} Institute of Mathematics. 

University of Natural Resources and Life Sciences, Vienna, Austria. 

{\sl E-mail:} {\tt gerald.kuba@boku.ac.at}

\end